# Data-Driven Surrogate Models for Agromaritime Applications: Finite Element–Neural Network Integration

## Muhammad Ilyas[1*]

[1] *Kuwait College of Science and Technology, Kuwait*
*Corresponding author's e-mail:* * **m.ilyas@kcst.edu.kw**

**Abstract**

Predicting nutrient transport and salinity distribution is crucial for mitigating climate-related threats to agromaritime systems. Traditional PDE-based models can capture the physics of nutrient dispersion, salinity and water quality. However, they face challenges in scalability and adaptability to real-time problems. In this article, we develop a hybrid approach that combines finite element discretisations with neural network integration to enable efficient and adaptive data-informed predictions. We use a finite element solver for the steady-state diffusion-reaction equation to generate a dataset across varying diffusivity, reaction and inflow conditions. We then build a proper orthogonal decomposition (POD), which reduces dimensionality, and a neural network (NN) that maps parameters to reduced coefficients. A numerical study presented on a simplified model demonstrates the proof-of-concept for predicting nutrient transport and salinity distribution. Numerical experiments show that the NN surrogate achieve a speed-up of approximately $956\times$ compared to a regular FEM solver while maintaining an accuracy of mean relative $L^2$-errors of $15\%$ across the test set, with occasional higher deviations, which is sufficient for rapid scenario screening and parametric studies. These results highlight the method's potential as a fast and accurate surrogate for nutrient and salinity prediction, offering a balance between FEM reliability and NN adaptability for sustainable agromaritime management.

**Keywords**: *Agromaritime Systems, Finite Element Method, Neural Networks surrogate modelling, Nutrient Transport, Proper Orthogonal Decomposition (POD)*.

## 1. Introduction

Sustainable agromaritime systems face significant challenges related to the management of water quality, salinity transport, nutrient dispersion, and associated flow dynamics. Accurate modelling of these processes is critical to ensure resilient agricultural and coastal practices, particularly in regions where freshwater and marine environments interact.

Salinity and nutrient transport are governed by complex physical processes. Numerical methods, such as the finite element method (FEM), have been widely applied to solve the underlying partial differential equations (PDEs) with high accuracy. Stabilised and mixed FEM formulations have been extensively studied in our previous works, including three-field Poisson equation (Ilyas & Lamichhane, 2016), sixth order problem (Droniou et al., 2019) as well as applications to complex physical systems such as wave-ice interactions (Ilyas et al., 2018). However, physics-based simulations remain computationally demanding, especially for parametric studies, scenario exploration, or real-time decision-making support.

92

In order to overcome these computational barriers, recent advances in artificial intelligence (AI) and data-driven modelling offer new opportunities to accelerate these simulations. Reduced-order models based on proper orthogonal decomposition (POD) have been successfully applied to nonlinear convection–diffusion systems (Fresca & Manzoni, 2022) and nonlinear structural analysis problems (Guo & Hesthaven, 2018). Physics-informed neural networks (PINNs) directly embed PDE structure into the learning process and have been proposed for transport in porous media and salinity modelling (Roh et al., 2023; Tartakovsky et al., 2020). Operator-learning methods such as DeepONet (Lu et al., 2021; Wang et al., 2021) and Fourier Neural Operators (FNO) (Hao & Song, 2024; Li et al., 2020) have emerged as powerful frameworks for parametric PDEs. Recent applications have also shown the promise of combining machine learning and reduced order modeling for PDEs (Tannous et al., 2025). These approaches have been extended to environmental systems including 1D advection–dispersion–reaction equations (Nguyen et al., 2025), 3D temperature–salinity anomaly reconstruction(Chen et al., 2023), seawater intrusion management using POD-based surrogates(Geranmehr et al., 2025), and PINN-based groundwater modeling (Ali et al., 2024).

Together, these works illustrate how AI-driven surrogates can accelerate PDE-based modelling of physical processes such as flow, transport, and reaction equation, making them particularly attractive for agromaritime applications where repeated simulation under varying environmental conditions is often required.

The primary contribution of this work is the development and assessment of a hybrid FEM-NN surrogate framework for agromaritime applications. Specifically, we integrate finite element simulations with a neural network (NN) and use Proper Orthogonal Decomposition (POD) to extract dominant solution features. This combination provides the reliability of FEM with the adaptability of NN to substantially reduce computational cost while maintaining acceptable accuracy for the fast prediction of parameterised PDEs. We demonstrate this approach on the steady-state diffusion–reaction equation, a representative PDE for salinity and nutrient transport. While this model neglects large-scale advection, it remains an important building block for understanding nutrient and salinity dynamics in controlled or microscale environments such as soil–root interactions, biofilms, and membranes. For larger-scale estuarine or aquifer systems, where advective transport dominates, the methodology can be extended naturally to the full ADRE in future work.

Compared to purely data-driven surrogates such as PINNs or DeepONets, the proposed POD–NN approach requires fewer data and computational resources, as it applies the physical structure embedded in the FEM-generated basis. While operator-learning models can achieve higher accuracy given extensive data, they typically require more training samples. The proposed hybrid framework, in contrast, balances interpretability, data efficiency, and computational cost, making it more suitable for engineering applications where limited but high-fidelity simulation data are available. The remainder of this paper is organised as follows. Section 2 presents the methodology, combining finite element discretisation, proper orthogonal decomposition, and neural network surrogates into an integrated research framework. Section 3 reports the numerical results, including representative snapshots, POD modes, and surrogate performance. Section 4 provides a discussion of the findings

93

and concludes with perspectives on the applicability of this hybrid approach to sustainable agromaritime systems.

## 2. RESEARCH METHODS

We will develop a quick prediction model for a salinity transport and nutrient dispersion by integrating machine learning techniques with diffusion-reaction numerical simulations. The research framework is illustrated in **Figure 1**.

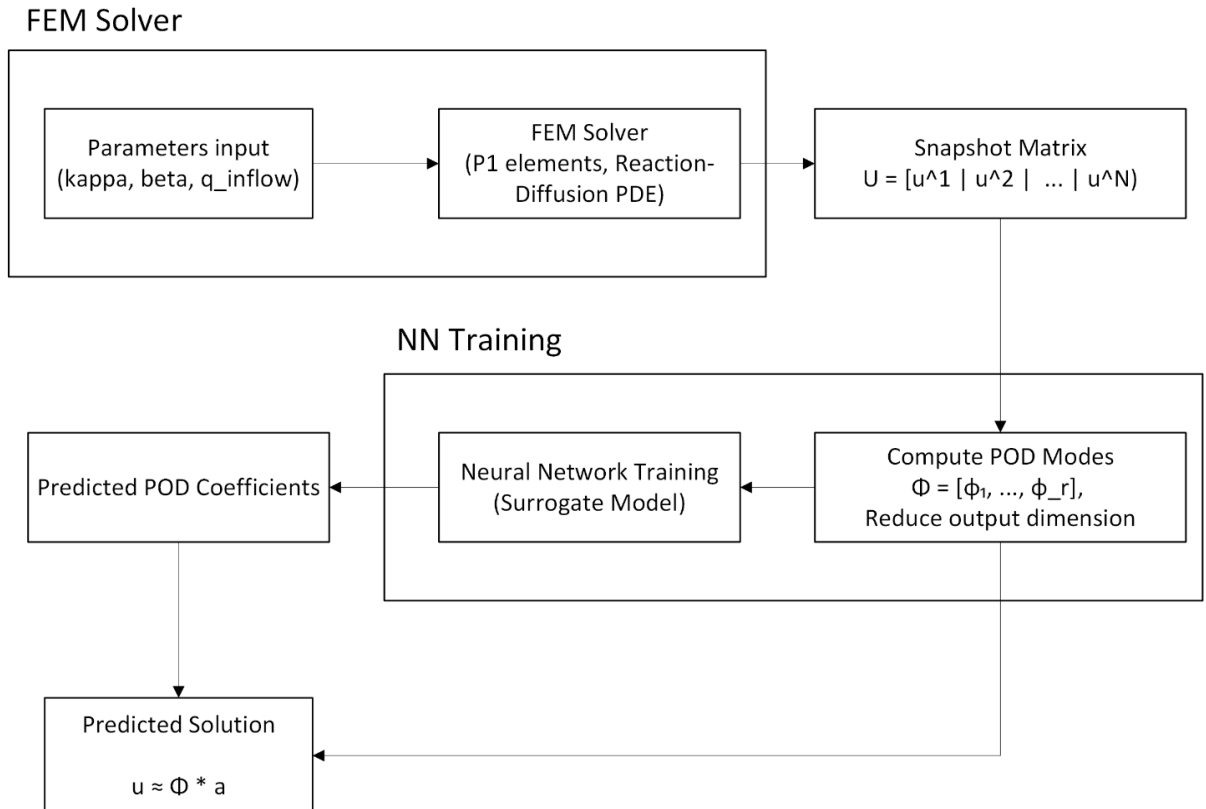

Figure 2: Workflow for the Hybrid FEM-NN Surrogate Model Construction

### 2.1 Governing Equation and Weak Formulation

We consider the steady-state diffusion–reaction equation

$$-\nabla \cdot (\kappa \nabla u) + \beta u = f \quad \text{in } \Omega, \tag{1}$$

where $u(x, y)$ is the solute concentration (e.g., salinity or nutrients), $\kappa > 0$ is the diffusion coefficient, $\beta \geq 0$ is the reaction coefficient, and $f$ is a source term. Boundary conditions are prescribed as Dirichlet values: $u = q_{in}$ on the left inflow boundary ($x = 0$), and $u = 0$ on all other boundaries.

The corresponding weak formulation is obtained by multiplying equation (1) by a test function $v$ in the Sobolev space $H^1(\Omega)$, integrating over $\Omega$, and applying the divergence theorem:

94

$$\kappa \int_\Omega \nabla u \cdot \nabla v \, d\Omega + \beta \int_\Omega u\, v, d\Omega = \int_\Omega f\, v, d\Omega + \kappa \int_\Gamma (\nabla u \cdot \boldsymbol{n}) v \, d\Gamma, \quad \forall v \in H_0^1(\Omega). \tag{2}$$

Here, $\boldsymbol{n}$ denotes the outward unit normal vector on the boundary $\Gamma = \partial\Omega$. The test function space is defined as $H_0^1(\Omega) = \{v \in H^1(\Omega) : v|_\Gamma = 0\}$, which simplifies the weak form by eliminating the boundary integral term when the test functions vanish on the Dirichlet boundary.

## 2.2 Finite Element Discretisation

Let $\mathcal{T}^h$ be a quasi-uniform triangulation of the polygonal domain $\Omega$ with a mesh size $h$. We use the standard linear finite element space $V_h \subset H^1(\Omega)$ defined on the triangulation $\mathcal{T}^h$

$$V_h = \{\eta \in C^0(\Omega) : \eta|_T \in \mathcal{P}_1(T), T \in \mathcal{T}^h\},$$

where $\mathcal{P}_1(T)$ denotes the space of linear polynomial on $T$.

Let $\{\phi_1, \phi_2, \ldots, \phi_N\}$ be a finite element basis function for $V_h$, which are piecewise linear function and satisfy the property $\phi_i(x_j) = \delta_{ij}$ for nodes $x_j$. We can write the approximate solution $u_h(x, y)$ as a linear combination of basis functions of $V_h$,

$$u_h(x, y) = \sum_{i=1}^N u_i \phi_i$$

where $u = (u_1, \ldots, u_N)^\top$ are to be determined.

Thus, equations (2) written in matrix form becomes:

$$(\kappa S + \beta M)u = f \tag{3}$$

where $S_{ij} = \int_\Omega \nabla\phi_i \cdot \nabla\phi_j \, dx$ is the stiffness matrix, $M_{ij} = \int_\Omega \phi_i \phi_j \, dx$ is the mass matrix, and $f_i = \int_\Omega f\, \phi_i \, dx$ is the load vector. The boundary conditions are incorporated into the system accordingly.

In order to generate the training dataset, physical and boundary parameters (diffusivity $\kappa$, reaction $\beta$, inflow strength $q_{in}$) are sampled across determined ranges. Each parameter triple $\mu^{(i)} = \left(\kappa^{(i)}, \beta^{(i)}, q_{in}^{(i)}\right)$ produces one FEM solution $u_h$ and the collection of all solutions gives the snapshot matrix:

$$U = [\boldsymbol{u}^{(1)} \quad \boldsymbol{u}^{(2)} \quad \ldots \quad \boldsymbol{u}^{(n_s)}] \in R^{N \times n_s}, \tag{4}$$

where $N$ is the number of degrees of freedom and $n_s$ is the number of snapshots.

## 2.3 Proper Orthogonal Decomposition

The snapshot matrix $\mathbf{U}$ is high-dimensional. We will apply proper orthogonal decomposition (POD) to reduce the dimensionality of the solution and select the dominant modes that captures the main dynamics of the solution manifold (Holmes et al., 2012; Quarteroni et al., 2016). Computationally, POD is equivalent to performing a singular value decomposition (SVD) of the snapshot matrix:

95

$$U = \Psi \Sigma V^\top, \tag{5}$$

where the columns of $\Psi = [\psi_1, \ldots, \psi_r] \in R^{N \times r}$ are the left singular vectors (the POD modes), $\Sigma = \text{diag}(\sigma_1, \ldots, \sigma_r)$ contains the singular values in descending order and $V$ contains the right singular vectors.

We will retain only the leading $m \ll n_s$ modes, chosen such that the retained singular values capture the dynamics of the solution based on the relative energy criterion:

$$\frac{\sum_{i=1}^{m} \sigma_i^2}{\sum_{i=1}^{r} \sigma_i^2} \geq \eta, \tag{6}$$

where $\eta$ is a chosen tolerance. The high threshold ensures the accuracy of the reduced basis while still achieving significant dimensional reduction. Using POD, any high-fidelity solution can then be approximated in this reduced basis:

$$u(\mu) \approx \sum_{j=1}^{m} a_j(\mu) \cdot \psi_j, \tag{7}$$

where $\mu = [\kappa, \beta, q_{in}]$ is the parameter vector and $a(\mu) = (a_1, \ldots, a_m)^\top \in R^m$ are the corresponding POD coefficients.

### 2.4. Neural Network Surrogate

The reduced-order representation obtained from POD is then mapped to parameters through a neural network (NN):

$$\mathcal{NN} : \mu \mapsto \hat{a}(\mu), \tag{8}$$

where $\hat{a}(\mu)$ is the prediction of the true POD coefficients $a(\mu)$.

Neural networks are universal function approximators (Goodfellow et al., 2016), and their recent use in physics-informed contexts has shown great potential for solving PDE-related problems (Raissi et al., 2019). By training the NN on the FEM-generated dataset, we obtain a surrogate that can predict new solutions for unseen parameter combinations at a fraction of the computational cost.

We implement the mapping (8) using a feedforward neural network. The network consists of 2 hidden layers with 100 neurons each. Hidden layers used hyperbolic tangent sigmoid activation functions, while the output layer used a linear activation. The network was trained using the Levenberg- Marquardt optimisation algorithm on 80% of the dataset and validated on the remaining samples generated from the FEM simulations and POD reductions. The learning process minimises the mean squared error between the predicted and true coefficients. The trained network can tphen predict the coefficients for any new parameter $\boldsymbol{\mu}^*$. The approximate solution is then recovered using equation (7):

$$u_{\text{NN}}(\mu^*) = \Psi_m \cdot \hat{a}(\mu^*), \tag{9}$$

where $\Psi_m = [\psi_1, \ldots, \psi_m]$ is the matrix of the first $m$ POD modes.

### 2.5 Integration Pipeline

The complete integrated framework is summarised as follows:
1. **High-Fidelity Simulation:** Generate the snapshot matrix $U$ by solving the FEM system (3) for a wide range of parameters.
2. **Dimensionality Reduction:** Compute the POD basis $\Psi_m$ and project all snapshots to obtain the coefficient dataset $\{(\mu^{(i)}, a^{(i)})\}$.
3. **Surrogate Training:** Train a neural network $\mathcal{NN}$ to approximate the mapping $\mu \mapsto a$.
4. **Fast Prediction:** For a new parameter $\mu^*$, evaluate the NN to get $\hat{a}(\mu^*)$ and reconstruct the approximate solution $u_{\text{NN}}(\mu^*)$.
5. **Evaluation:** The performance evaluated by comparing the approximate solution $u_{\text{NN}}(\mu^*)$ with FEM reference solution $\boldsymbol{u_h}$ using relative $L^2$-errors of representative fields.

This integration of FEM, POD, and NN provides a lightweight surrogate model that accelerates PDE simulations while preserving essential physical structures for parameterised PDEs.

### 3. Results and Discussion

In this section, we will illustrate the numerical framework showing several snapshots along the way. We define the problem on $\Omega = [0,10] \times [0,5]$ with gridsize $200 \times 100$. We have the Dirichlet boundary condition as: $u = q_{\text{in}}$ on the left inflow boundary ($x = 0$), and $u = 0$ on all other boundaries. We also set the source term as $f = \sin(\pi x) \cdot \sin(\pi y)$. We generate 500 sample solutions for the training and trial data with the parameters randomly selected across predetermined ranges: $\kappa \in [10^{-3}, 10^{-1}], \beta \in [0,1], q_{in} \in [0.1,1]$. Some sample FEM solution snapshots are shown in Figure 2 illustrating the variation in solution patterns.

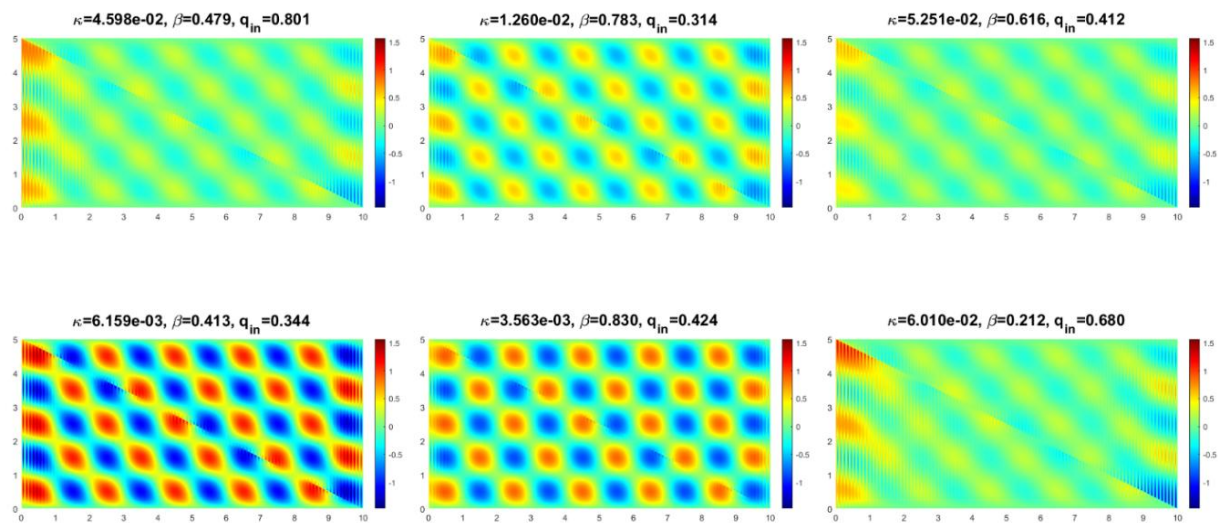

Figure 2: Sample FEM Solution Snapshots

Using the FEM solution produced, we calculate a POD basis and project all snapshots

97

to obtain the coefficient dataset. In this study, we set that the selected modes capture at least $99.9\%$ of the energy contained within the snapshot matrix. The cumulative sum of the first 6 POD modes are shown in Figure 3. Based on this result, we can see that we can select just the first 4 modes to capture $99.9\%$ of the main dynamics of the solution. The dominant 4 modes are shown in Figure 4. Using this criterion, we employed a neural network described in Section 2.4 on the POD coefficient dataset. In order to evaluate the training result, we compared the POD solution with the NN reconstruction. Figure 5 shows the histogram distribution of the relative $L^2$-error with mean and median around $16\%$ and the worst case deviations above $61\%$. We also show the best and the worst case NN reconstruction in Figure 6.

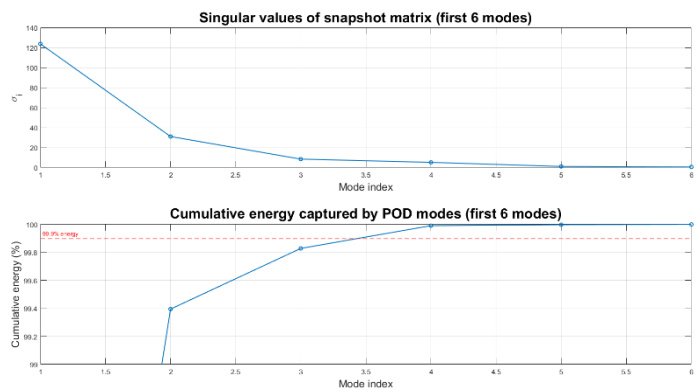

Figure 3. Singular Value of Snapshot Matrix and Cumulative Energy Sum

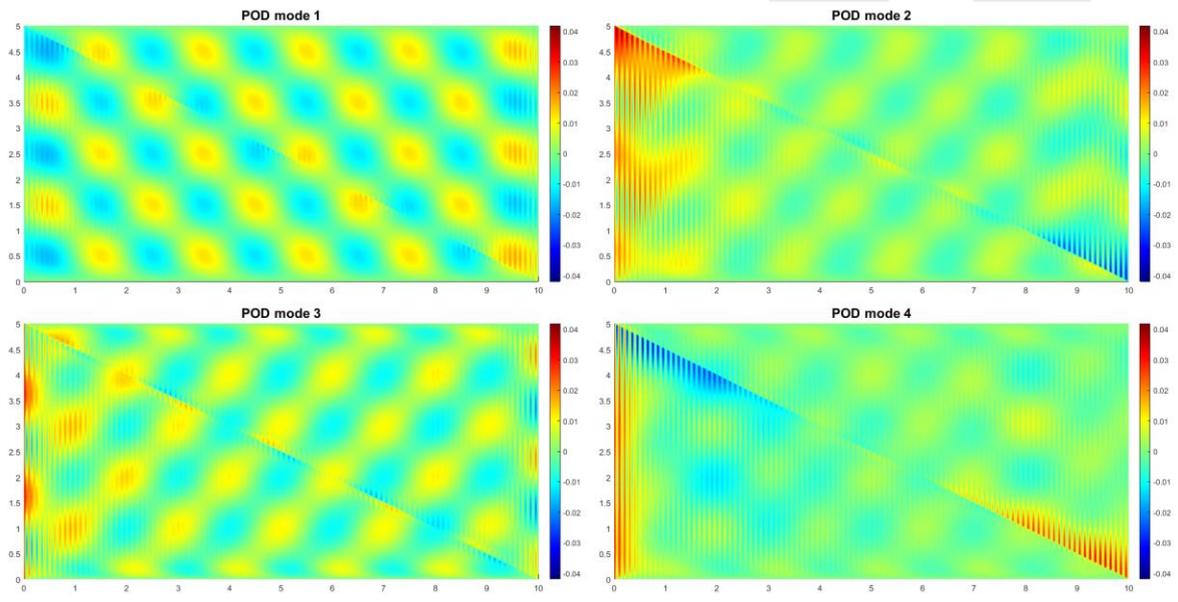

Figure 4. Leading 4 POD Modes

98

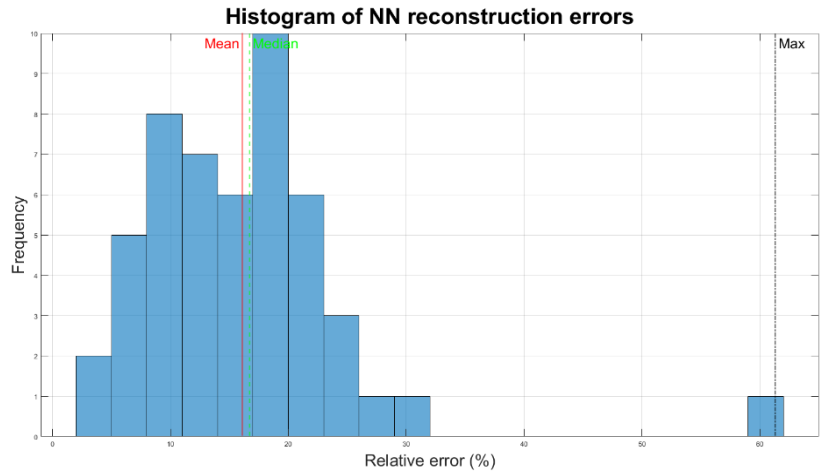

Figure 5. Histogram of Neural Network (vs POD) reconstruction error

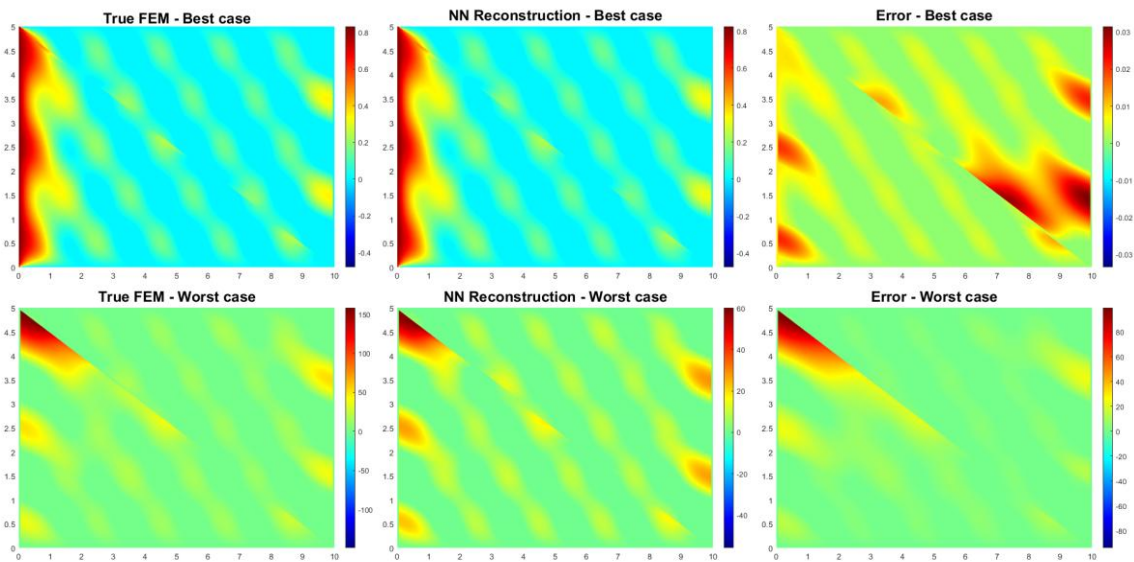

Figure 6. Best ($\kappa = 8.171 \times 10^{-02}, \beta = 0.815, q_{in} = 0.833$) and Worst ($\kappa = 1.101 \times 10^{-03}, \beta = 0.001, q_{in} = 0.100$) NN reconstruction

We believe that the mean relative $L^2$-error of approximately 16% arises from two primary sources: (i) the projection error introduced by POD basis truncation, and (ii) the regression error of the neural network approximation. The worst-case deviations exceeding 60% (**Figure 6**) typically occur for parameter combinations near the boundaries of the sampled space or in regions where the solution manifold exhibits sharp gradients that challenge both the POD reduction and NN approximation. Despite these errors, the accuracy remains sufficient for rapid parametric studies and scenario screening, where identifying trends and relative behaviours is often more critical than obtaining high-fidelity solutions.

As an evaluation benchmark, we compared the efficiency of the actual FEM solver with the NN reconstruction using random parameters within the specified ranges and domains described. The results are summarised in **Table 1**. We also illustrate a 3D sample of the NN reconstruction comparison in **Figure 7**. These results demonstrate

99

6 – 7 November 2025 | Hybrid Conference
IPB International Convention Center (IICC), Bogor, Indonesia & Online via Zoom

that the POD-NN framework retains physical structure while enabling fast parametric approximation.

Table 1. FEM vs NN runtime data

|  |  |
|---|---|
| Average FEM time | 5.850223 seconds |
| Average NN time | 0.006119 seconds |
| Speed-up factor (FEM/NN) | 956.1x |
| Mean relative error | 15.22 % |
| Median relative error | 15.33 % |
| Max relative error | 33.29 % |

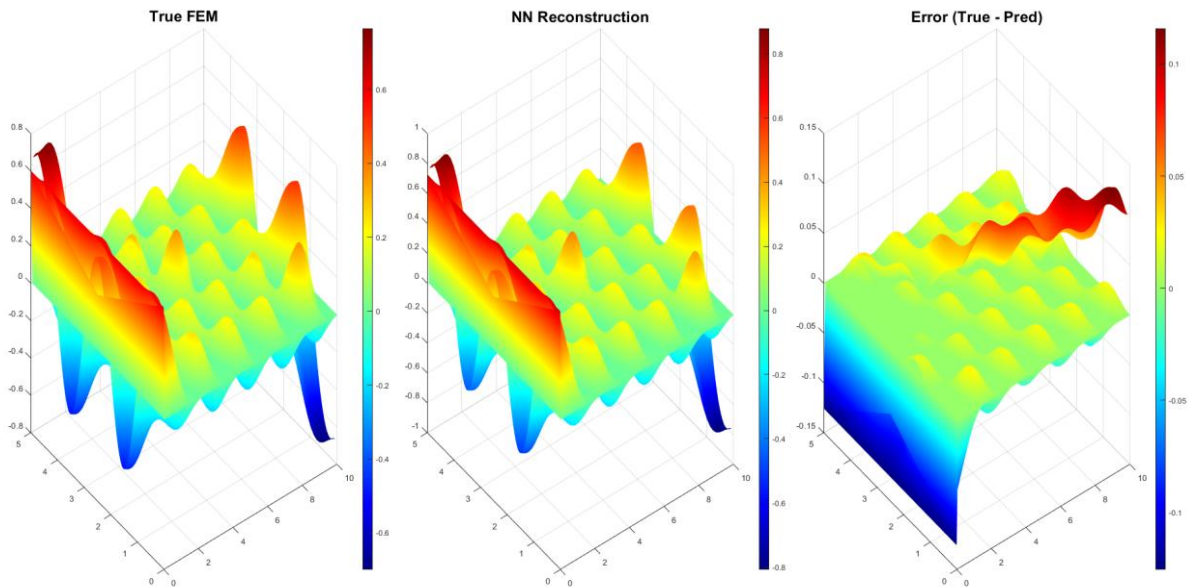

Figure 7. True FEM vs NN reconstruction sample ($\kappa = 5.639 \times 10^{-02}, \beta = 0.559, q_{in} = 0.603$)

**Remark:** *The histogram in Figure 5 shows a mean relative error of 16% when comparing the NN reconstruction against the POD projection. In the runtime comparison against the actual FEM solver (Table 1), the mean error was 15.22%, demonstrating consistency between the reduced-order model and high-fidelity solutions.*

## 4. Conclusion

We developed and evaluated a hybrid FEM-NN framework for nutrient and salinity transport modelling. By training a neural network on FEM-generated datasets and compressing the solution space with POD, we achieved surrogate predictions with mean relative errors of 15% while reducing computational time by approximately 956 times faster compared to direct FEM solvers. These results confirm that data-driven surrogates can provide rapid and reasonably accurate approximations suitable for parametric studies and near real-time applications.

Compared with fully data-driven operator-learning models such as PINNs, DeepONets, or FNOs, the proposed POD–NN framework is more data-efficient and

100

requires significantly less training while maintaining strong interpretability through its FEM-derived basis. This balance between physics fidelity and computational efficiency underlines the novelty of the approach and its practicality for engineering-scale agromaritime problems.

However, the main limitation is the simplified test case: a steady-state diffusion-reaction equation on a simple rectangular domain. While this study illustrates feasibility, future work should extend the framework to the full advection-dispersion-reaction systems and more realistic geometries. Extending the approach to spatially varying parameters or exploring operator-learning architectures such as DeepONet or FNO would further enhance its accuracy and robustness.

In conclusion, the proposed FEM–NN integration offers a promising pathway toward fast, adaptive, and physically consistent surrogates for agromaritime systems. By balancing physics-based rigour with data-driven efficiency, this approach supports scenario exploration and decision-making in environments where computational speed and resilience are equally critical.